\documentclass[A4, 12pt]{amsart}
\usepackage{enumerate}

\let \hx \hbox
\let \vx \vbox

\def\C{\mathbb{C}}
\def\E{\mathbb{E}}
\def\A{\mathbb{A}}

\def\F{\mathbb{F}}

\def\R{\mathbb{R}}
\def\Z{\mathbb{Z}}
\def\G{\mathbb{G}}

\def\hfl#1#2{\smash{\mathop{\hx to 12mm{\rightarrowfill}}
\limits^{\scriptstyle#1}_{\scriptstyle#2}}}
\def\vfl#1#2{\llap{$\scriptstyle #1$}\left\downarrow
\vx to 6mm{}\right.\rlap{$\scriptstyle #2$}}

\def\build#1_#2^#3{\mathrel{\mathop{\kern 0pt#1}\limits_{#2}^{#3}}}

\newtheorem{thm}{Theorem}[section]
\newtheorem{prop}{Proposition}[section]
\newtheorem{lem}{Lemma}[section]
\newtheorem{cor}{Corollary}[section]

\newtheorem{defi}{Definition}[section]

\theoremstyle{remark}
\newtheorem{rem}{Remark}

\def\g0{G^{(0)}}

\def\dx{\partial X}
\def\intm{\overset{\:\circ}{M}}
\def\intf#1{\overset{\:\circ}{#1}}
\def\gr{grou\-poid}
\def\d{{\rm d}}

\def\cstar{$C^*$-algebra}

\def\ie{{\it i.e.}}
\def\cinf{$C^\infty$}
\def\M{M}
\def\G{{\bf G}}
\def\E{\mathcal{E}}

\def\sv{submanifold}
\def\ll{longitudinally smooth}
\def\lb{\hbox{lb}}
\def\rb{\hbox{rb}}
\def\pd{pseudodifferential}
\def\var{manifold}
\def\tts{that's to say\ }
\def\equa{{decoupage}}
\def\mc{manifold with corners}

\def\S{\mathcal S}
\def\A{\mathcal A}

\begin{document}

\author{Bertrand Monthubert}
\title[Pseudodiff. calculus on  manifolds with corners and groupoids]{Pseudodifferential calculus on  manifolds \\with corners and groupoids}
\date{July 1997}
\address{Institut de Math\'ematiques\\Universit\'e Pierre et Marie Curie (Case 191)\\4, place  Jussieu\\F-75252 PARIS CEDEX 05 }
\email{monthube@math.jussieu.fr}

\subjclass{Primary 35S15, 47G30, 22A22; Secondary 46L80, 19K56}
\keywords{Pseudodifferential calculus, manifolds with corners,
  groupoids, $C^*$-algebras}

\renewcommand{\abstractname}{R\'esum\'e}
\begin{abstract}

     Nous construisons un groupo\"\i de diff\'erentiable longitudinalement
  lisse associ\'e \`a une vari\'et\'e \`a coins. Le calcul pseudo-diff\'erentiel
  sur ce groupo\"\i de co\"\i ncide avec le calcul pseudo-diff\'erentiel
  de Melrose (aussi appel\'e $b$-calculus). Nous d\'efinissons \'egalement
  une algèbre de fonctions \`a d\'ecroissance rapide sur ce groupo\"\i de
  ; il contient les noyaux des op\'erateurs r\'egularisants du (petit)
  $b$-calculus. 
 
\end{abstract}
\maketitle
\section{Introduction}
Let $X$ be a \cinf\ \var\ modelled on $\R_+^n$; for
any point $m \in X$ the number of zero coordinates is called the codimension of $m$.

An open face $\intf{F}$ is a connected component of a set of points of same
codimension. The closure $F$ of an  open face is called a face. 

Following R. Melrose (see \cite{melrose-aps}), $X$ is said to be a \mc\ if every hyperface is embedded, more
 precisely if for every hyperface $F$ there exists a definition function
 $\rho_F \geq 0$ on $X$, which vanishes on $F$ only, with  $\d \rho_F \neq 0$
on $F$, and such that in every intersection  of several
 hyperfaces the corresponding $\d \rho_F$ are independent.

In \cite{melrose-aps}, R. Melrose  introduced a \pd\ calculus on \var s with boundary, and more generally with corners.
Through the ``$b$-stretched product'' construction, he defined a space of \pd\ operators adapted to boundary problems, the $b$-calculus.

It appears in the work of Alain Connes (see  \cite{connes}) that a \pd\ calculus is naturally associated to any Lie \gr. The (usual) \pd\ 
calculus on a smooth \var\ without boundary  corresponds to the Lie \gr\ $M\times M$. Moreover, introducing Connes' tangent \gr\ 
helps to better understand Atiyah-Singer's index theorem (\cite{atiyah-singer}).  Pursuing these ideas, we have made this construction explicit 
for any differentiable \gr\ \ll\ in \cite{bm-fp}. 

The goal of this article is to show that Melrose's $b$-calculus is the \pd\ calculus on a suitable \gr. To a \var\ with corners $X$, we  
associate a  differentiable \gr\ $\Gamma(X)$, which is \ll, and we show that the \pd\ calculus on this \gr\ is precisely  the $b$-calculus.
We also give a definition of rapidly decreasing functions  on $G$ (with respect to a proper homomorphism of \gr s with values in  $\R^p$), 
which allows to consider additional smoothing operators.

Our construction generalizes straight forwardly to certain \gr s with corners (the \gr\ $\Gamma(X)$ above corresponds to the \gr\ $M 
\times M$ associated to a \var\  without boundary $M$).

\bigskip
We begin by introducing, in paragraph 2, some additional material
about \pd\ calculus on a differentiable \gr, the main ingredient being a Schwartz space on 
a \gr\ associated with a proper homomorphism with values in $\R^p$. 

In paragraph 3, we build the \gr\ associated to a smooth \var, with respect to a \sv\ of codimension 1, then we generalize to a transverse 
family of \sv s of codimension 1.

Finally, in paragraph 4 we apply the results of the previous sections to the context of a manifold with corners. We show that the $b$-calculus 
coincides with the \pd\ calculus on the \gr\ associated to the manifold.

I wish to thank Fran\c cois Pierrot, Michael Pimsner and Georges
Skandalis for helpful discussions.

\begin{rem}
  Some of
these results have been independently obtained by V. Nistor, A. Weinstein
and P. Xu (see \cite{nwx}).
\end{rem}

 \section{Pseudodiffential calculus on differentiable groupoids}

 \subsection{Schwartz kernels associated to \pd\ operators on a
   differentiable \gr}
Let $G$ be a differentiable \gr\ \ll, \tts a  \gr\ together with a structure of
\cinf\  manifold with corners on $G$, $\g0$ being a \cinf\ \sv\ (with corners) of $G$, such that the fibers $G^x=r^{-1}(x), x \in 
\g0$ are smooth \var s,  the  range  and source maps, $r$ and $s$, are submersions,  and the composition is \cinf.  

In \cite{bm-fp} we showed that for any differentiable \gr\ \ll, one can build a \pd\ calculus.

We are going to give an equivalent definition, in terms of Schwartz kernels, analogous to the definition of  H\"ormander
(\cite{hoermander}) in the case of \pd\  operators on a smooth \var.

Let's define the Lie algebroid associated to a differentiable \gr\  \ll, $G$, denoted by  $\A(G)$, as the normal bundle of $\g0$ in $G$,
 $N_\gamma\g0=T_\gamma G/T_\gamma \g0$. It is canonically isomorphic to the restriction to the set of units of $G$ of the kernel of the 
morphism $\d s \hbox{(or } \d r):TG \to T\g0$; the  sections of $\A(G)$ are the vector fields $v$ on $G$ which are left-invariant and 
such that 
 $\d s(v(\gamma))=0, \forall \gamma \in G$. 

The vector bundle of half-densities $\Omega^\frac{1}{2}$ is the line bundle  over $G$ whose fiber in $\gamma \in G$ 
is the linear space of maps
$$\rho : \Lambda^kT_\gamma G^{r(\gamma)} 
\otimes \Lambda^kT_\gamma G_{s(\gamma)} \to \C$$
such that $\rho(\lambda\nu)=|\lambda|^{1/2} \rho(\nu), \forall
\lambda \in \R$. Here, $k$ is the dimension of the fibers  $G^{r(\gamma)}$
and $G_{s(\gamma)}$. This bundle is trivial: a choice of a metric trivializes it.

Let $\varphi$ be a diffeomorphism from a (tubular) neighborhood $O$ of $\g0$ in
 $G$ to $\A(G)$, such that the following diagram commutes:
$$\begin{array}{ccc}
O & \hfl{\varphi}{} & \A(G) \cr
\vfl{r}{} & & \vfl{}{} \cr
\g0 & \hfl{}{} & \g0. 
\end{array}$$

We moreover assume that 
$\varphi(\gamma)=0 \iff \gamma \in \g0$,
and  that its differential is the isomorphim $T_xG^x\simeq N_x \g0$.

\begin{defi}
The space of  order $m$ \pd\ kernels on a \gr\ $G$ is the space
$I^m(G,\g0;\Omega^\frac{1}{2})$ of distributional sections $K$ on $G$ with values in
 $\Omega^\frac{1}{2}$ which are \cinf outside $\g0$, and given by the oscillatory integral in
 a neighborhood of $\g0$ inside $O$,
$$K =(2\pi)^{-n}\int_{\A_{r(\gamma)}(G)^*} 
e^{i\langle \varphi(\gamma),
 \xi \rangle} a(\gamma,\xi) \d  \xi,$$
where  $a$ is a polyhomogeneous symbol  of order $m$ with  values in $\Omega^\frac{1}{2}$. 

We denote by $I^\infty(G,\g0;\Omega^\frac{1}{2})$ the space of  \pd\ kernels   of any order 
($I^\infty(G,\g0;\Omega^\frac{1}{2})= \cup_{m\in \Z} I^m(G,\g0;\Omega^\frac{1}{2})$). 
\end{defi}
Let's remark that smooth sections of $\A(G)$ are in $I^1(G,\g0;\Omega^\frac{1}{2})$.
Of course one could also choose more general symbols, like symbols of type $(\rho,\delta)$.

There is a (convolution) pairing $I^\infty \times C_c^\infty \to C^\infty$. $I^\infty$ is not an algebra, but 
$I_c^\infty(G,\g0;\Omega^\frac{1}{2})$, the space of compactly supported kernels of \pd\ operators, is an involutive algebra. 
Moreover, the convolution of an element in $I_c^\infty $ and an element of $I^\infty $ is an element of $I^\infty $.

\bigskip 
The definition above is equivalent to that of \cite{bm-fp},
definition 2.2, where we considered the \pd\ operators on $G$ as $G$-equivariant families of operators
on the fibers $G_u$. 
Indeed, a $G$-operator $T$ of order $m$ is given by a family of \pd\ operators of order $m$ $T_u$, equivariant with respect to 
the natural action of $G$. Each $T_u$ is given by its Schwartz kernel $K_u \in I^m(G_u^2, \Delta_u, 
\Omega_u^{\frac{1}{2}})$, thus the equivariance relation implies: 
$$K_{s(\gamma)}(\gamma ',\gamma'' )=K_{r(\gamma)}(\gamma'\gamma^{-1},\gamma''\gamma^{-1})$$
for $\gamma',\gamma'' \in G_{s(\gamma)}$.

\bigskip 
Such a family allows to define an element of  $I^m(G,\g0;\Omega^\frac{1}{2})$ by  
$K(\widetilde \gamma)=K_{s(\gamma)}(\gamma,\gamma')$ where $\gamma$ and $\gamma'$ are such that  
$\gamma'\gamma^{-1}=\widetilde \gamma$. The definition is independent of the choice of $\gamma$ and
$\gamma'$ : let $\delta, \delta'$ be such that $\delta'\delta^{-1}=\widetilde \gamma$. Thus let 
$\alpha=\gamma^{-1}\delta$. The equivariance relation gives:
$$K_{s(\delta)}(\delta,\delta')=K_{s(\alpha)}(\delta,\delta')=K_{r(\alpha)}(\delta\alpha^{-1}, 
\delta'\alpha^{-1})=K_{s(\gamma)}(\gamma,\gamma').$$

Conversely, one can define a family of operators from an element of
\linebreak $I^m(G,\g0;\Omega^\frac{1}{2})$ by  
$K_{s(\gamma)}(\gamma,\gamma')=K(\gamma'\gamma^{-1}).$
\bigskip

\begin{rem}
  It seems that the definitions in \cite{nwx} are equivalent to the
  ones above. They have been obtained independently.
\end{rem}

 \subsection{Rapidly decreasing functions}

We have seen that Schwartz kernels of \pd\ operators are smooth outside the diagonal. Thus, outside a neighborhood of the diagonal, 
we just work with \cinf\ functions.

To obtain an algebra, we ask {\it a priori} that the functions be compactly supported. However, if we are given a proper \cinf\ map, 
$\varphi : G \to \R^p$, which is in addition a morphism of \gr s (in the sense that
$\varphi(\gamma_1\gamma_2)=\varphi(\gamma_1)+\varphi(\gamma_2)$), we can then define a notion of rapidly decreasing 
functions on $G$.
Let's remark that such a function does not necessarily exist; for
instance, its existence immediately implies that 
$\g0 \subset \varphi^{-1}(0)$ is compact.

If we are given such a morphism $\varphi$, we can define the space 
$$\S^0=\{ f \in C^{\infty}(G) , \forall P \in \C[X_1,\ldots,X_p],
 \sup_{\gamma \in G} \|P(\varphi(\gamma)) f(\gamma)\| < \infty \}.$$

\begin{defi}
The space $\S(G)$ of rapidly decreasing functions on $G$, or
\linebreak Schwartz space on $G$, is the subspace of $\S^0$ 
of functions such that $$\forall \kappa\in I_c^\infty(G,\g0;\Omega^\frac{1}{2}), \kappa * f \in \S^0, f * \kappa \in 
\S^0.$$

\end{defi}
One can easily check that $\S(G)$ is a subalgebra of $C^*(G)$.

Let's note that to define $\S(G)$ it is actually enough to consider
the action of differential operators.

Thus we can define, in this context, a \pd\ calculus where functions are taken in
 $\S(G)$, not only in $C_c^{\infty}(G)$.

\begin{rem}
This definition coincides with the usual one in the case where $G$ is $\R^n$, considered as a Lie group (whence a \gr\ with $\g0$ a point).
\end{rem}

The kernels of operators that we consider are sums of kernels of operators with compact support and of rapidly decreasing 
functions on $G$; they form an algebra $I_s^\infty(G,\g0;\Omega^\frac{1}{2})$ for the convolution product.

\section{The groupoid associated to a transverse family  of codimension 1 submanifolds}

\subsection{The puff of a  codimension 1 submanifold}

Let $\M$ be a \var,  and $V\subset \M$ be a \sv\  of $\M$ of codimension 1. We can define 
the normal bundle to $V$ in $\M$, that we denote by $NV$, as the quotient $TM/TV$ (for $x \in V, N_xV=T_xM/T_xV$).

We  denote by $\G(V)$ the \gr\ composed of the disjoint union of $\M\backslash V \times
\M\backslash V$ and of the fiber bundle $\widetilde{V \!\!\times\!\! V}$ on $V \times V$ of the isomorphisms
 from $N_{y}V$ to  $N_{x}V$, glued together in the following way: 

$$(\M\backslash V)^2 \ni (x_n,y_n) \to (x,y,\alpha) \in\widetilde{V \!\!\times \!\!V} \!\!\iff\!\! 
\left\{ \begin{array}{l}
   \!\! x_n\to x, y_n\to y, \cr
   \!\! p_{x}(x_n-x)=\alpha
    (p_{y}(y_n-y))\cr
    + o(p_y(y_n-y)) \cr
\end{array} \right. $$
where $p_{x}$ is the projection of $T_{x}\M$ on the factor $N_{x}V$.

The   \gr\ structure on $\M\backslash V \times \M\backslash V$ is that of the \gr\ of the
trivial equivalence relation, and on
$\widetilde{V\!\!\times\!\! V}$
the source and range maps, and the composition are defined as:
$$ \left\{\begin{array}{l}
        r(x,y,\alpha)=x\cr
        s(x,y,\alpha)=y\cr
        (x,y,\alpha) (y,z,\alpha')=(x,z,\alpha\circ\alpha').
        \end{array} \right.$$ 

The differentiable manifold structure is given by the following charts. On
$M\backslash V\times M\backslash V$ the structure is given by that of $M$.

Now, for $i=1,2$, let $\varphi_i$ be a diffeomorphism from an open $\Omega_i$
of $M$  to $(V \cap \Omega_i)\times \R$, such that if $x \in V \cap
\Omega_i,\varphi_i(x)=(x,0)$.  The $\varphi_i$ determine a trivialisation of $NV$ over $V\cap \Omega_i$.

Then $W_{\varphi_1,\varphi_2}=(V \cap \Omega_1)\times (V \cap
  \Omega_2)\times \R \times \R^*$ defines a chart of $\G(V)$ by
  $\psi(x,y,t,\lambda)= \varphi^{-1}_1(y,t)\times
\varphi^{-1}_2(x,\lambda t)$ if $t\neq 0$, and by 
$\psi(x,y,0,\lambda)= (x,y,\alpha)$ where $\alpha$ is the isomorphism from
  $N_yV$ to $N_xV$ which, in the trivialisations given by the $\varphi_i$, is the  dilation of rate $\lambda$.

It is clear that $\G(V)$ is a Lie \gr.

\bigskip
In addition, the \gr\ structure being given by  $r \oplus s(x,y,t,\alpha)=(x,y,t,\alpha t)$, we have $d(r \oplus s) 
(x,y,t,\alpha)(\xi,\eta,\tau,\iota)=(\xi,\eta, \tau, \alpha \tau+ t\iota)$ and in $t=0$ (\tts  on $V^2$), 
$d(r \oplus s) (x,y,0,\alpha)(\xi,\eta,\tau,\iota)=(\xi,\eta, \tau, \alpha \tau)$. Thus 
$\forall \gamma \in \widetilde{V \!\!\times\!\! V}, 
d(r \oplus s)(T_\gamma   \G(V)$ is the hyperplan of
 $T_{r \oplus s(\gamma  )} {\g0}^2$ defined by the equation 
$\pi_1(\xi)=\alpha(\pi_2(\xi))$, where $\pi_2$ is the projection of the second coordinate of $\xi$ on $N_{s(\gamma)}(V)$, and 
$\pi_1$ is the projection   of the first coordinate of $\xi$ sur $N_{r(\gamma)}(V)$.

\bigskip
Two particular situations can be stressed right now:
\begin{enumerate}[a)]
\item {\bf if $V$ is transversally oriented}
Let's remark that in this case the normal bundle can be trivialised, the bundle 
 $\widetilde{V \!\!\times\!\! V}$ being trivialised in  $V \!\!\times\!\! V \times \R^*$.

Moreover, the normal bundle can be decomposed into a positive part and a negative part. Thus it is possible to restrict the bundle 
$\widetilde{V \!\!\times\!\! V}$ to the isomorphisms which preserve each subbundle. We denote by $\G(V)^+$ the \gr\  
obtained; it is open in $\G(V)$. By using a tubular neighborhood of $V$ in $M$, we can describe $\G(V)$ with two charts:
 $\M\backslash V \times \M\backslash V$ and $V\times V\times \R \times \R^*$.

\item {\bf if $V$ splits $\M$ in two parts}

In this context, we can restrict  $\G(V)$ by keeping only $M^+\!\!\times\!\! M^+ \cup M^-\!\!\times\!\!M^-$, glued with the 
bundle of the isomorphisms in $\widetilde{V \!\!\times\!\! V}$ which preserve each subbundle.
We denote by $\G(V_+)$ the \gr\  obtained; it is open and closed in $\G(V)$.

If we provide a trivialisation of  $N(V)$, $\G(V_+)$ is the union of 
$\M_+\backslash V \!\!\times\!\! \M_+\backslash V \cup \M_-\backslash V \!\!\times\!\! \M_-\backslash V$ 
and $V \!\!\times\!\! V \times \R^*_+$.

\end{enumerate}

\subsection{Puff of a transverse family of  codimension 1 submanifolds}

If we consider several  codimension 1 \sv s of $\M$, we can make a similar construction provided certain conditions are fulfilled.
Let's first of all recall a few facts about transversality and fiber products.

Let $X_1, \ldots, X_n, Y$ be a  family of \cinf\ \var s. A family of \cinf morphisms 
   $\varphi_i : X_i \to Y$ is called transverse if, for any family $(x_i)_{1 \leq i \leq n} \in  \Pi_{1 \leq i \leq n}
   X_i$
   such that $\exists y \in Y, \forall i, \varphi_i(x_i)=y$, the
   orthogonals (in $T_y^*Y$) of the spaces $\d \varphi_i(T_{x_i}X_i)$
   are in direct sum.

Under this condition, the fiber product of the $X_i$ over $Y$ is a \cinf\ submanifold of  $\Pi_{1 \leq i \leq n}   X_i$.

Besides, a family of \sv s will be called transerse if the inclusions of these \sv s are transverse.

\begin{defi}
We call \equa\ a \var\  $\M$ together with a finite family 
 $(V_i)_{i \in I}$ of codimension 1 \sv s such that $\forall J \subset I$, the family of inclusions of $V_j, j\in J$ is transverse.
\end{defi}

From now on we denote $V_J=\cap_{i\in J} V_i$, for $J\subset I$.

An important consequence of this definition is the following proposition:
\begin{prop}\label{equa-transverse}
Let $(\M,(V_i)_{i \in I})$ be a \equa. 
Then the differentiable morphisms $r \oplus s : \G(V_i) \to \M^2$ are transverse.
\end{prop}

An immediate corollary of the proposition above is that the fiber product of the  $G(V_i)$ over $\M^2$ is a \cinf\ \var.

\begin{defi}
 Let $\E=(\M,(V_i)_{i \in I})$ be a \equa. The puff of the \equa\ 
 in $\M^2$, denoted by $\G(\E)$, is the fiber product of the  morphisms $r \oplus s: \G(V_i) \to \M^2$. 
\end{defi}
It is a  sub\gr\ and a \sv\ of $\Pi_{i \in I} \G(V_i)$, thus it is a Lie \gr.

It is clear to that point that every \sv\ $V_J$ is saturated in $\G(\E)$, \tts that every element of $\G(\E)$ whose source is 
in $V_J$ has its range in $V_J$.
The orbits of  $\G(\E)$ are the
$$\intf{F_J}=\left( \bigcap_{i\in J} V_i \right) \cap 
\left( \bigcap_{i\not\in J} M\backslash V_i \right).$$

\subsection{Puff in a transverse groupoid}
We are now going to generalize the process by defining the puff of a \gr\ along a transverse \equa.

\begin{defi}
Let $G$ be a Lie \gr; set $M=\g0$.   We say that $G$ is transverse to
a \equa\ $(\M,(V_i)_{i\in I})$ if for any   $J \subset I$,  the    
differentiable morphisms
$r\oplus s:G \to \M^2$ and the inclusions  $V_i^2$ in $\M^2$, for $i\in J$, are transverse.
\end{defi}

The following easy proposition gives some conditions  equivalent to this definition:
\begin{prop}
  Let $\E=(\M,(V_i)_{i\in I})$ be a \equa, and $G$ be a Lie \gr\  such that $\g0= \M$ .
Then the following  are equivalent:
\begin{enumerate}[(i)]
\item $G$ is transverse to  $\E$;
\item for every $x \in \M$, and for every $J \subset I$,  the differentiable morphisms  $s : G^x \to
\M$ and the inclusion of $V_J$ in $\M$ are transverse;
\item the morphisms $r \oplus s : G \to \M^2$ and $r \oplus s: \G(\E) \to \M^2$ are transverse.
\end{enumerate}
\end{prop}

It is thus possible to ``puff'' a \equa\ in $G$:

\begin{defi}
We define the puff of the \equa\ in  $G$, 
denoted by $\langle G : \E \rangle$, as the fiber product:
$$\begin{array}{ccc}
\langle G : (V_i) \rangle & \hfl{\phi}{} & G \cr
\vfl{}{} & & \vfl{}{r \oplus s} \cr
\G(\E) & \hfl{r \oplus s}{} & \M \times \M. 
\end{array}$$

\end{defi}

We get immediately the following proposition, due to the associativity of the fiber product: 

\begin{prop}
Let $G$ be Lie \gr, transverse to a \equa\ \linebreak[4] $\E=(\M,(V_1,\ldots,V_n))$. Then
$$\langle G : \E\rangle = 
\langle \ldots \langle\langle G : V_1 \rangle : V_2 \rangle \ldots \rangle.$$
In particular,
 $$\G(\E)=\langle \ldots \langle\langle M\!\!\times\!\!M :
 V_1 \rangle : V_2 \rangle \ldots \rangle.$$
\end{prop}
We can thus obtain $\G(\E)$ by puffing successively the \sv s $V_i$ from $M\times M$.

\begin{cor}
  For every permutation $\sigma$, 
$$\langle \ldots \langle\langle G : V_1 \rangle : V_2 \rangle \ldots
\rangle = \langle \ldots \langle\langle G : V_{\sigma_1} \rangle : V_{\sigma_2} \rangle \ldots
\rangle.$$

\end{cor}
We can thus puff the \sv s following any order.

\bigskip
Let's examine what happens for a single codimension 1 \sv.

Let $V$ be a codimension 1 \sv\ of a \var\ $\M$, and $G$ be a Lie \gr\ transverse to
$V$. Then the fiber product $G'$ of $G$ and $\G(V)$ over $\M^2$ equals the union of $G_{M \backslash V}^{M \backslash V}$ 
and of the bundle $\widetilde{G_V^V}$
of the isomorphisms of $s^*NV$ to  $r^*NV$ over $G_V^V$, both \gr s being glued together as for $\G(V)$.

Indeed if $\gamma' \in G'$, then $\gamma'=(\gamma,\gamma_V)$, with $\gamma \in G, \gamma_V \in \G(V)$ such that
 $r \oplus s (\gamma)=r \oplus s (\gamma_V)$. Then necessarily, 
$r \oplus s (\gamma) \in V^2$, \tts  $\gamma \in G_V^V$, or $r \oplus s (\gamma) \in (M \backslash V)^2$, \tts 
$\gamma \in G_{M \backslash V}^{M \backslash V}$. In the first case, $\gamma'$ is given by  an isomorphism of $N_{y}V$ to  $N_{x}V$, \tts $\gamma' \in \widetilde{G_V^V}$. In the second case, $\gamma'$  
is completely determined by $\gamma \in G_{M  \backslash V}^{M \backslash V}$, 
so $\gamma' \in G_{M \backslash V}^{M \backslash V}$.

\bigskip
To conclude this section, let's remark that the construction of the puff of a \equa\ in a  transverse \gr\ is canonical. It does not depend, in 
particular, of a  trivialisation of the normal bundle  (in the case of a  transversally oriented \sv). However, if we are provided with such a  
trivialisation, we can simplify the form of the  \gr. Let's consider the case where the \sv s are transversally oriented.

In this case, each normal bundle $N(V_i)$ can be trivialised, thanks
to a definition function of $V_i$, denoted by $\rho_i$. 
Then the bundle $\widetilde{V_i \!\!\times\!\! V_i}$  trivialises in $V_i \!\!\times\!\! V_i \times \R^*$, and 
$$G(\E)=\bigcup _{J \subset \{1,\ldots,k\}} \bigcap_{j \in J}V_j^2 \times (\R^*)^{card(J)}. $$
The topology on $G(\E)$ is given by the convergence of the following sequences:
 $$\bigcap_{j \in J}V_j^2 \times (\R^*)^{card(J)} \ni (x_n,y_n,\lambda_n) 
\hfl{}{n \to +\infty}  (x,y,\mu) \in \bigcap_{j' \in J'}V_{j'}^2 \times (\R^*)^{card(J')},$$ where $J'=J \cup 
\{j_0\}$. There is convergence if and only if:
$$\left\{ \begin{array}{l}
x_n \to x \\
y_n \to y \\
\lambda_j \to \mu_j, j \in J \\
\frac{\rho_{j_0}(y_n)}{\rho_{j_0}(x_n)} \to \mu_{j_0}.
\end{array} \right. $$

We obtain a \gr: the source and range maps, and the 
composition are defined like this:
$$ \left\{\begin{array}{l}
        r(x,y,\lambda)=x\cr
        s(x,y,\lambda)=y\cr
        (x,y,\lambda)\circ (y,z,\lambda')=(x,z,\lambda\lambda')
        \end{array} \right.$$ 
if $x,y$ and $z$ are in the same maximal intersection of \sv s.
If we wish to restrict to the isomorphisms which preserve each sub-bundle of the normal bundle, we obtain the same result by replacing each time 
$\R^*$ by $\R^*_+$.

Finally, we see that the puff allows to get, from a \gr\ transverse to a family of \sv s, a \gr\ which is longitudinal 
on these \sv s. This is especially interesting in the context of manifolds with corners, which comes out of the case when each
 $V_i$ separates $\M$ in two.

\subsection{The groupoid associated to a manifold with corners}

Every \var\ with corners has a smooth neighborhood, thus we will take
 the following definition, which is equivalent to that of the introduction:

\begin{defi}
A \mc\ $X$ is the restriction of an oriented \equa, $(\M,(V_i)_{i \in
 I})$, where each $V_i$ separates $\M$ in two, to the positive part of
 $\M$.
The faces are the restrictions of the intersections of the \sv s. We
 say that $X$ is the positive part of this oriented \equa.
\end{defi}

\bigskip
From now on, $X$ is a compact \mc.

We showed in \cite{bm-fp} how to construct a \pd\ calculus far any
\ll\ differentiable \gr. In the case of a \mc, we can not build  a \pd\ calculus
on the \gr\ $X \times X$, because for every $x \in X, 
(X\!\!\times\!\! X)_x=X$ which is not smooth.

However, this \gr\ can be desingularized by considering the positive
part of the puff of the \equa\ $\E=(M,(V_i)_{i \in I})$ for which $X$
is the positive part, following an approach due to Melrose.

Indeed, as each $V_i$ separates $M$ in two, we can build the fiber
product of the  $\G({V_i}_+)$, that we denote by
$\G(\E_+)$. 

Let's note that the positive part of the \equa\
is saturated in $\G^+(\E)$: indeed, $\Gamma(X)=(\G(\E)^+)_X^X$ is such
that if $x \in \intf{F}, \Gamma(X)^x=\intf{F}\times
(\R^*_+)^{codim(F)}$ which is a \cinf\ \var. Thus it is \ll.
Also, this \gr\ is independent from the \equa\ considered: if $X$
embeds in another \equa, the positive part remains unchanged, and
$\Gamma(X)$ as well for it only depends on the positive part.

\begin{rem}
The \cstar\ of $\Gamma(X)$ is of type I and its irreducible
representations are indexed by the open faces of $X$, which correspond
to subquotients isomorphic to  $C_0(\R^k,\mathcal{K}_{\intf{F}})$,
where $\mathcal{K}_{\intf{F}}$ is the
\cstar\ of the compacts of $L^2(\intf{F})$. Thus it is easy to describe the spectral sequence which gives the $K$-theory of 
$C^*(\Gamma(X))$. It coincides with that of \cite{Melrose-Nistor}, and we get simple proofs of these results.
\end{rem}
\begin{rem}
  Nistor, Weinstein and Xu defined a  \gr, in the case of a manifold
  with boundary,  in \cite{nwx}, which seems to be isomorphic
  to $\Gamma(X)$. It was obtained independently.
\end{rem}

\section{Pseudodifferential calculus on manifolds with corners}

\subsection{Schwartz space on manifolds with corners}

We are now provided with all the material to build a \pd\ calculus with
compact support on manifolds with corners. 

To be able to consider, in addition to compactly supported \pd\ operators, some ``rapidly decreasing'' smoothing operators, we need to define a 
proper morphism of \gr s, $\Gamma(X) \to \R^k$, where $k$ is the number of  \sv s of the \equa\ (\tts the maximal  codimension of the 
faces of $X$). 

Let's first of all consider a smooth manifold $M$ together with a \sv\ $V$ which splits $M$ in two. Thanks to
the definition functions of V, $\rho_V$, we can define $\psi_V : \G(V_+) \to \R$ by: 

\begin{itemize}
\item if $\gamma \in \intm\times\intm, \psi_V(\gamma)=\log(\frac{\rho_V(r(\gamma))}{\rho_V(s(\gamma))})$;

\item otherwise, $\gamma=(x,y,\lambda)$ and $\varphi_V(\gamma)=\log(\lambda)$.
\end{itemize}

It is easy to check that $\varphi_V$ is a proper map. Now, for a \mc\ $X$, we can consider a \equa\ $(M,(V_i)_{i\in I})$ 
whose positive part is $X$. The morphisms $\varphi_{V_i}$ then induce a proper morphism of \gr s at the level of the fiber product: 
$\varphi:\Gamma(X) \to \R^{card(I)}$.

This allows us to consider a Schwartz space on $\Gamma(X)$.
\begin{prop}
If we change the definition functions of the faces, the difference
between the new morphism $\varphi'$ and $\varphi$ is a bounded
map. Thus it does not change the space of rapidly decreasing functions.
\end{prop}
\begin{proof}
Let's denote by $\rho'_i$ the new definition functions of the
faces. Then there exists \cinf\ functions $f_i$ such that $\rho'_i=f_i \rho_i$, with $f_i>0$. So $\d\rho'_i(x)=f_i(x)\d\rho_i(x)$, and 
every
isomorphism of $N_yV_i$ to $N_xV_i$ which associates $\d\rho'_i(x)$
to $\d\rho'_i(y)$ equals $f_i(x)/f_i(y)$ times the isomorphism which  associates $\d\rho_i(x)$ to $\d\rho_i(y)$. Finally, 
$$\varphi'_i(\gamma)=\log(f_i(r(\gamma))/f_i(s(\gamma))) + \varphi_i(\gamma).$$
But $f_i$ is continuous and does not vanish, and is defined on a
compact, so the  difference $\varphi'_i - \varphi_i$ is bounded.
\end{proof}
This allows to build a canonical \pd\ calculus:
\begin{defi}
Let $X$ be a \mc. The (canonical) \pd\ calculus on $X$ is the \pd\
calculus (with rapidly decreasing functions) on the \gr\ $\Gamma(X),$
with respect to $\varphi$.
\end{defi}

\begin{rem}[Pseudodifferential calculus on differentiable \gr s with corners]
We can easily generalize this method to \gr s such that:
\begin{enumerate}
\item  $G$ and $\g0$ are manifolds with corners
\item there exists a differentiable \gr\ (without corners) $G'$ which contains
  $G$, such that ${G'}^{(0)}$ admits a structure of \equa, $G'$ is
  transverse to this \equa, and ${G'}^{\g0}_{\g0}=G$.

\end{enumerate}
We shall call {\it differentiable \gr\ with corners} such a \gr.

In this context we can define a \pd\ calculus by considering the puff
in $G'$ of the \equa\ defined on ${G'}^{(0)}$, then by taking its
 restriction to $\g0$, this \gr\ being \ll. 
\end{rem}

\subsection{Relation with the $b$-calculus}

For a connected \var\ with boundary, the \gr\ 
$$\Gamma(X)=\intf{X} \times \intf{X} \cup \bigcup_{1 \leq i \leq l} \partial_i X\times\partial_i 
X\times\R^*_+,$$
 where $i$ indexes the boundary components, leads to the following
 result:

\begin{thm}
The \pd\ calculus on $\Gamma(X)$ coincides with the \linebreak (small) $b$-calculus of Melrose.
\end{thm}
\begin{proof}
Let's recall shortly the main points of $b$-calculus
(\cite{melrose-aps}).  
To simplify the notation, let's assume that the boundary is
connected.

The classical \pd\ calculus can be defined in terms of Schwartz
kernels, by considering the \pd\ operators as distributional sections
on $X^2$ with values in the bundle of half-densities, which are \cinf\
on $X^2 \backslash \Delta$. In the $b$-calculus, the space $X^2$ is
replaced by the $b$-stretched product, $X_b^2$ which is defined as the union of $X^2 \backslash (\dx)^2$ and of the
hemispheric normal bundle of  $(\dx)^2$, \tts the 
quotient of the normal bundle of $(\dx)^2$  by the natural
action of $\R^*_+$. The \cinf\ structure is such that if we fix a definition function of the boundary, $\rho$, the map 
$\tau(x,y)=\frac{\rho(x)-\rho(y)}{\rho(x)+\rho(y)}$, defined on $\intf{X}^2$, extends to a \cinf\ map on $X_b^2$.
We also have a natural map $\beta_b : X_b^2 \to X^2$. 

This new space has some boundary components, in particular 
$$\lb=\hbox{cl}(\beta_b^{-1}(\dx \times \intf{X}))=
\tau^{-1}(-1), \rb=\hbox{cl}(\beta_b^{-1}(\intf{X} \times
\dx))=\tau^{-1}(1)$$
 (cl stands for the closure). 
Finally, the $b$-\pd\ operators on a \var\ with boundary $X$ are
defined by their  Schwartz kernels, as  distributional sections smooth outside the diagonal
on the $b$-stretched product, which vanish, as well as all their
derivatives, on $\lb \cup \rb$.

\bigskip
Now let's remark:
\begin{lem}
Let $X$ be a \var\ with boundary, and $\Gamma(X)$ be the \gr\ associated to
$X$. Then $\Gamma(X)$ is an open \sv\ of $X_b^2$, and 
$$X_b^2\backslash \Gamma(X) = \lb \cup \rb.$$ 
\end{lem}
\begin{proof}

On $\Gamma(X)$, we  actually have $\tau(\gamma)=\tanh
(\varphi(\gamma))$, so $\Gamma(X)= \linebreak \tau^{-1}(]-1,1[)$; moreover 
$\tau$ is a submersion so $\Gamma(X)$ is an open \sv\ of $X_b^2$, of complementary $lb \cup rb$. 
\end{proof}

We can deduce from this lemma that, as $lb \cup rb$ is disjoint from a neighborhood of $\Delta_b$, the \pd\ calculus with compact support 
on $\Gamma(X)$ coincides with the  distributional sections  on $X_b^2$, smooth outside the diagonal
on the $b$-stretched product $X_b^2$, which vanish on a neighborhood of $lb \cup rb$.
\end{proof}
\begin{rem}
  The fact that compactly supported \pd\ operators are $b$-\pd\
  operators was  announced independently in \cite{bm-fp} and in  \cite{nwx}.
\end{rem}
To completely describe the (small) $b$-calculus of Melrose it remains to understand the nature of the smoothing operators, \ie\ the \cinf\ 
functions on $X_b^2$ which vanish in Taylor series on $lb\cup rb$.

Let  $f$ be a \cinf\  function on a compact neighborhood of a point of $lb\cup rb$
adapted to $\rho$, and $g$ be its restriction to $\Gamma(X)$.  If $f$ vanishes
with all its derivatives on $lb\cup rb$,  Taylor formula implies that
$g$ is rapidly decreasing with respect to the powers of $\tau(x,y) \pm
1$, which is equivalent to $\exp(\pm \varphi)$. These are thus {\it a
  fortiori} rapidly decreasing functions with respect to $\varphi$. This means that $f \in \S^0$. Now, as the small $b$-calculus is an algebra, 
if $P$ is a (compactly supported) \pd\ operator on $\Gamma(X)$ (thus a $b$-\pd\ operator), $Pf $ is a smoothing operator of the $b$-
calculus, and $Pf\in \S^0$. Thus $f \in \S$. So the smoothing operators of the $b$-calculus are in $\S(G)$. More precisely, they are the 
\cinf\ function which are rapidly decreasing with respect to $\exp(\varphi)$.

\begin{rem}
In the case of a manifold with corners, we can build a space
  $X_b^2$ for each \sv\ of a \equa\ whose positive part is $X$. Then we can make the fiber product of these spaces, and finally consider 
the restriction to the positive part. It gives a $b$-stretched product for the \mc.
\end{rem}

\bibliographystyle{amsplain}
\bibliography{tg}

\end{document}